\title{\large Computation of Cohomology Operations\\ on Finite Simplicial Complexes
\thanks{
Authors are partially supported by the  PAICYT research project FQM-296 from
Junta de Andalucia and the DGES--SEUID research
project
PB98--1621--C02--02
from Education and Science Ministry (Spain).}}
\author{\normalsize R. Gonz\'{a}lez--D\'{\i}az, P. Real\\
\normalsize Dept. of Applied Math.,\\
\normalsize University of Seville,  Spain\\
 \{rogodi, real\}@us.es
 }
 \date{}

\documentstyle[emlines2,bezier,12pt]{article}
\newcommand{\ra}{\rightarrow}

\newcommand{\D}{\displaystyle}
\newcommand{\scst}{\scriptscriptstyle}
\newcommand{\scsz}{\scriptsize}

\newcommand{\ot}{\otimes}

\newcommand{\ti}{\times}

\begin{document}

\maketitle

\begin{abstract}
We propose a method for calculating cohomology operations
for finite simplicial complexes.

Of course, there exist well--known methods for computing (co)homology groups, for
example,  the ``reduction algorithm" consisting in reducing the
matrices  corresponding to the differential in each dimension
to the Smith normal form, from which one can read
off (co)homology groups of the complex \cite{Mun84}, or the  ``incremental
algorithm" for computing Betti numbers \cite{DE93}. However,
there is a gap in the literature concerning general methods
for computing cohomology operations.

For a given finite simplicial complex $K$, we sketch a procedure including the computation of
some primary and
secondary cohomology operations and
the $A_{\infty}$--algebra structure on the cohomology of $K$.
This method is based on the transcription of the reduction
algorithm mentioned above, in
terms of a special  type of algebraic homotopy equivalences,
called a contraction, of the (co)chain complex of $K$ to a
``minimal" (co)chain complex $M(K)$. For instance, whenever
 the ground ring is a field or the (co)homology of $K$ is free,  then $M(K)$ is isomorphic to
  the (co)homology
of $K$.
Combining  this contraction with
the combinatorial formulae for Steenrod reduced $p$th powers  at cochain level
developed in \cite{GR99} and \cite{Gon00},  these
operations at cohomology level can be computed.
Finally, a method
for calculating Adem secondary cohomology operations
$\Phi_q:\; Ker(Sq^2H^q(K))\ra H^{q+3}(K)/Sq^2H^q(K)$ is showed.
\end{abstract}

\section{Introduction}

Particular important  topological invariants are the (co)homology
groups. In a certain way, these groups measure the degree of
connectedness of the space. Although there are plenty of programs
for calculating (co)homology groups of finite simplicial
complexes, we have not found any general software for computing
cohomology operations.

Our  main motivation is the design of a program for computing all
sort of cohomology invariants on finite simplicial complexes:
(co)homology groups, cup product, Bockstein cohomology operation,
cohomology operations determined by homomorphisms of coefficient
groups, Steenrod squares and reduced $p$th powers, Pontrjagin
squares and $p$th powers, the $A_{\infty}$--algebra structure of
cohomology, higher cohomology operations, etc. In this paper, we
give a solution to the problem of computing Steenrod squares and
reduced $p$th powers \cite{Ste47,ES62} and Adem secondary
cohomology operations \cite{Ade52,Ade58}. Our approach is based
on the translation of the well--known ``reduction" algorithm for
computing (co)homology groups \cite{Mun84} in terms of homotopy
equivalences. In that way, we have a description of the
generators of the (co)homology groups in terms of cochains. In
fact, this is sufficient to enable us to determine the effect of
the induced maps between cohomology groups corresponding to
cochain maps.
Using the same approach
we think that the rest of primary cohomology operations could be
attacked.

\section{Background} \label{prel}

We give a brief summary of concepts and notation used in the following sections.
Our terminology follows Munkres \cite{Mun84}.

For $0\leq q\leq n$, a $q$--{\em simplex} $\sigma$ in ${\bf R}^n$ is the
convex hull of a set $T$ of $q+1$
affinely independent points $(v_0,...,v_q)$. The dimension of $\sigma$ is $|\sigma|=q$.
For every non--empty $U\subset T$, the simplex $\tau$ defined by $U$ is a  {\em face} of $\sigma$.
A  {\em simplicial complex} $K$ is a collection of simplices   satisfying the following
 properties:

\begin{itemize}
\item If $\tau$ is a face of $\sigma$ and $\sigma\in K$ then $\tau\in K$.
\item If $\sigma, \tau\in K$ then $\sigma\cap \tau$ is either empty or a face of both.
\end{itemize}
The set of all the $q$--simplices of $K$ is denoted by $K^{(q)}$.
The largest dimension of any simplex in $K$ is the {\em dimension of $K$}.
A simplex $\sigma$ in $K$ is {\em maximal} if it is not face of any simplex in $K$.
Therefore, $K$ can be given by the set of its maximal simplices.
A subset $L\subset K$ is a {\em subcomplex} of $K$ if it is a simplicial complex itself.
All simplices in this paper
have finite dimension and all simplicial complexes are finite collections.
From now on, $K$ denotes a finite simplicial complex.
The {\em oriented} $q$--simplex $\sigma=[v_0,...,v_q]$ is the equivalence class of the particular
ordering $(v_0,...,v_q)$. Two orderings are equivalent if they differ from one another by an even
permutation.

Let $\Lambda$ denote an abelian group.
A formal sum,  $\lambda_1 \sigma_1+\cdots+\lambda_n\sigma_n$,
where $\lambda_i\in \Lambda$ and  $\sigma_i$ are oriented $q$--simplices, is called a $q$--{\em chain}.
The {\em chain complex} canonically associated to $K$, denoted by  $C_*(K)$, is the family of groups
such that in each dimension $q$, $C_q(K)$ is the group of $q$--chains in $K$.
The {\em boundary} of a $q$--simplex $\sigma=[v_0,v_1,...,v_q]$ is the $(q-1)$--chain

$$\partial_q \sigma=\sum_{i=0}^q(-1)^i[v_0,v_1,\dots,\hat{v}_i,\dots, v_q]\,,$$
where the hat means that $v_i$ is omitted. By linearity, the boundary operator $\partial_q$
can be extended to $q$--chains, where it is a homomorphism.
 It is clear that for each $q$--simplex $\sigma_j$ there exist unique integers $\lambda_{ij}$
such that

$$\D\partial_q(\sigma_j)=\sum_{\tau_i\in K^{(q-1)}}\lambda_{ij} \tau_i\,.$$
The matrix $A_q=(\lambda_{ij})$ is the {\em matrix of} $\partial_q$ relative to the
bases $K^{(q)}$ and $K^{(q-1)}$.
The group of $q$--{\em cycles}, $Z_q(K)$, is the kernel of $\partial_q$, and define
$Z_0(K)=C_0(K)$. The group of
$q$--{\em boundaries},
$B_q(K)$, is the image of $\partial_{q+1}$, that is, the subgroup of $q$--chains $b\in C_q(K)$ for
which there exists a $(q+1)$--chain $a$ with $b=\partial_{q+1} a$.
It can be shown that $\partial_q\partial_{q+1}$ is
null so $B_q(K)$ is a subgroup of $Z_q(K)$. Then, the $q$th {\em homology group}

$$H_q(K)=Z_q(K)/B_q(K)$$
can be defined for each integer $q$.
Let $K$ and $L$ be two simplicial complex.
A {\em chain map} $f: C_*(K)\rightarrow C_*(L)$ is a family of homomorphisms

$$\{f_q: C_q(K)\rightarrow C_q(L)\}_{q\geq 0}$$
such that $\partial_qf_q=f_{q-1}\partial_q$, for all $q$.

Dual concepts to the previous ones can be defined.
 The {\em cochain complex} canonically
associated to $K$, $C^*(K)$, is the family

$$C^*(K)=\{C^q(K),\delta_q\}_{q\geq 0},$$
where

$$C^q(K)=\mbox{Hom}(C_q(K);\Lambda)=\{c: C_q(K)\rightarrow \Lambda,\quad\mbox{ $c$ is a homomorphism}\}$$
and $\delta^q:C^q(K)\ra C^{q+1}(K)$ called the  {\em coboundary operator} is given by

$$\delta^q (c)(a)=c(\partial_{q+1}a)\,,$$
where $c\in C^q(K)$ and $a\in C_{q+1}(K)$. Observe that a $q$--cochain can be
defined only on $K^{(q)}$ and extended to $C_q(K)$ by linearity.
Moreover, if $\Lambda$ is a ring, then
a basis of $C^q(K)$ is the set of homomorphisms

$$\sigma^*:C_q(K)\rightarrow \Lambda\,,$$
such that if $\tau\in K^{(q)}$, then $\sigma^*(\tau)=1$ if $\tau=\sigma$, and $\sigma^*(\tau)=0$
otherwise.
$Z^q(K)$ and $B^q(K)$ are the kernel of $\delta^q$ and the image of $\delta^{q-1}$, respectively.
The elements in $Z^q(K)$ are called $q$--{\em cocycles} and those in $B^q(K)$
are called $q$--{\em coboundaries}. It is also satisfied that $\delta^q\delta^{q-1}=0$ so
the $q$th {\em cohomology group}

$$H^q(K)=Z^q(K)/B^q(K)$$
can also be defined for each integer $q$.
If $\Lambda$ is a ring,
the cohomology of $K$ is also a ring with the {\em cup product}

$$\smile: H^p(K)\ot H^q(K)\ra H^{p+q}(K)$$
defined at cocycle level by

$$c\smile c'(v_0, v_1,\dots,v_{p+q})=\mu (c(v_0,\dots,v_p)\ot c'(v_p,\dots,v_{p+q}))\,,$$
where
$v_0<v_1<\cdots<v_{p+q}$, $c$ is  an $p$--cocycle, $c'$ is a $q$--cocycle
and $\mu$ is the product on $\Lambda$.

We use in this paper a special type of homotopy equivalences.   A
{\em contraction} $r$ of a  chain complex $N_*$ to another
chain complex $M_*$ is a set of three homomorphisms $(f,g,\phi)$
where $f:N_n\ra M_n$ ({\em projection}) and $g:M_n\ra N_n$ ({\em inclusion})
are chain maps
and satisfy that $fg=1_{\scst M}$, and
 $\phi: N_n\ra N_{n+1}$ ({\em homotopy operator})  satisfies that

$$1_{\scst N}-gf=\phi \partial_{\scst N}+\partial_{\scst N}\phi\,.$$
Moreover, $\phi g=0\,,\quad f\phi=0\,,\quad \phi\phi=0\,.$
A contraction {\em up to dimension} $n$ of $N_*$ to $M_*$
consists in a set of three homomorphisms  $(f,g,\phi)$ such that

$$f_k:N_k\ra M_k\,,\qquad g_k:M_k\ra N_k\qquad\mbox{ and }\qquad \phi_{k-1}: N_{k-1}\ra N_{k}$$
are defined for all $k\leq n$,  $\phi_n=0$, and the conditions of being a contraction
are satisfied up to dimension $n$.
Starting from a contraction $r=(f,g,\phi)$ of $N_*$ to $M_*$,
 it is possible to give another contraction
$r^*=(f^*,g^*,\phi^*)$  of
Hom$(N;\Lambda)$ to Hom$(M;\Lambda)$ as follows:

$$f^*:\mbox{Hom}(N_n;\Lambda)\ra\mbox{Hom}(M_n;\Lambda)\,,\qquad
g^*:\mbox{Hom}(M_n;\Lambda)\ra\mbox{Hom}(N_n;\Lambda)\,,$$

$$\phi^*:\mbox{Hom}(N_n;\Lambda)\ra\mbox{Hom}(N_{n-1};\Lambda)\,,$$
are such that

$$f^*(c)=cg\,,\qquad g^*(c')=cf\qquad\mbox{ and }\qquad \phi^*(c)=c\phi\,,$$
where $c\in \mbox{Hom}(N_n;\Lambda)$ and  $c'\in \mbox{Hom}(M_n;\Lambda)$.

\section{``Minimal" Chain Complexes}

It is possible to translate the results of the
``reduction algorithm", discussed at length in \cite{Mun84}, in terms of homotopy equivalences.
Combining this translation with  modern homological perturbation techniques,
 algorithms for computing
algebraic invariants, such as the $A_{\infty}$--algebra structure on the cohomology of $K$
and primary and
secondary cohomology operations can be designed
in an easy way.

First of all, it is necessary to recall the reduction algorithm for computing homology groups
of a finite simplicial complex $K$. This method
consists in
reducing the matrix $A$ of  the boundary operator in each dimension  $q$, relative to given bases
 of $C_q(K)$ and  $C_{q-1}(K)$,
to its Smith normal form $A'$ (a matrix of integers
satisfying that all its elements  are zero except for
$\lambda'_{11}\geq 1$ and $\lambda'_{11}/\lambda'_{22}/\cdots/\lambda'_{\ell\ell}$
for some integer $\ell$). This reduction is done in each dimension $q$
modifying the given base of $C_{q-1}(K)$, using the following
``elementary row operations'' on the matrix $A$:

\begin{enumerate}
\item[(1)] Exchange row $i$ by row $k$.
\item[(2)] Multiply row $i$ by $-1$.
\item[(3)] Replace row $i$ by row $i+n($row $k)$, where $n$ is an integer and $k\neq i$.
\end{enumerate}
Of course, there are similar ``column operations" on $A$ corresponding to changes of basis of
$C_{q}(K)$.
With this operations, the Smith normal form $A'$ of $A$ can be obtained,
relative to some bases
$\{a_1,\dots a_r\}$ of $C_q(K)$ and $\{e_1,\dots, e_s\}$
of $C_{q-1}(K)$.
Then,

\begin{enumerate}
\item[(1)] $\{a_{\ell+1},\dots,a_r\}$ is a basis of $Z_q(K)$,
\item[(2)] $\{\lambda'_{11}e_1,\dots,\lambda'_{\ell\ell}e_{\ell}\}$ is a basis of $B_{q-1}(K)$.
\end{enumerate}

Obviously, a dual treatment for $C^*(K)$ and, consequently, for the cohomology $H^*(K)$,
can be done.

A chain complex $M_*(K)$ is called {\em minimal} if
in each dimension $q$, $M_q(K)$ is a finitely generated
free abelian group and
the Smith normal form $A'$ of the differential of $M_q(K)$
 has the first element $\lambda'_{11}$  different from $1$.
  An {\em algebraic minimal model} of $K$
is a minimal chain complex $M_*(K)$ together with a contraction of $C_*(K)$ to
$M_*(K)$. Indeed, there is an algebraic minimal model for any finite simplicial complex $K$
and any two algebraic minimal models of $K$ are isomorphic.

Now, let us construct inductively an algebraic minimal model of a given finite simplicial complex $K$.
Suppose that an algebraic minimal model up to dimension $q-1$
is already constructed.
That is, we have a minimal chain complex  $M'_*(K)$ such that $M'_i(K)=0$, $i\geq q$, and a contraction
up to dimension $q-1$,
$(f',g',\phi')$,
of $C_*(K)$ to $M'_*(K)$.
Reduce the matrix of $\partial_q:C_q(K)\ra C_{q-1}(K)$ to its Smith normal form $A'$. If
the elements $\lambda'_{11}=\cdots=\lambda'_{tt}=1$, for $t\leq \ell$ (that is, $\partial(a_i)=e_i$
for $1\leq i\leq t$), then
define $M_*(K)$ as follows:

$$\begin{array}{l}
M_i(K)=M'_i(K),\qquad \mbox{for } i\neq q-1,q\\\\
M_{q-1}(K)=M'_{q-1}(K)-\Lambda[e_1,\dots,e_t]\\\\
M_q(K)=C_q(K)-\Lambda[a_1,\dots,a_t]
\end{array}$$
where $\Lambda[a_1,\dots,a_t]$ and
$\Lambda[e_1,\dots,e_t]$ are the free abelian groups generated by $\{a_1,\dots,a_t\}$
and $\{e_1,\dots,e_t\}$, respectively.
The formulae for the component morphisms of the contraction up to dimension $q$, $(f,g,\phi)$,
of $C_*(K)$ to $M_*(K)$ are:

\begin{eqnarray*}
&&f(x)=\left\{\begin{array}{ll}
f'(x)\qquad&\mbox{if } x\in \Lambda[e_{t+1},\dots,e_s]\;\mbox{ or }\;x\in C_i(K), \;i<q,\\
0 &\mbox{if }x\in \Lambda[e_1,\dots,e_t]\;\mbox{ or }\; x\in \Lambda[a_1,\dots,a_t],\\
x &\mbox{if }x\in \Lambda[a_{t+1},\dots,a_r],
\end{array}\right.\\\\
&&g(y)=\left\{\begin{array}{ll}
g'(y)\qquad&\mbox{if } y\in M_i(K), \;i<q,\\
y&\mbox{if } y\in M_n(K),\end{array}\right.\\\\
&&\phi(x)=\phi'(x)\qquad\mbox{if } x\in C_i(K), \;i<q-1,\\
&&\phi(e_i)=a_i \qquad\mbox{if } 1\leq i\leq t,\\
&&\phi(e_i)=0 \qquad\mbox{if } t+1\leq i\leq s.
\end{eqnarray*}

In this way, we can determine an algebraic minimal model for a finite simplicial complex $K$.
Observe that whenever $\Lambda$ is a field or the homology of $K$ is free, then $M_*(K)$
is isomorphic to  $H_*(K)$ and, therefore, we can obtain
a contraction of $C_*(K)$ to its homology.

Passing to cohomology does not represent a problem and a dual process can be done without effort.

The fact of dealing with contractions is highly important
in obtaining topology invariants such as the $A_{\infty}$--algebra structure
of the cohomology of $K$ \cite{GS86}. In particular, if $\Lambda={\bf Q}$, then
 from the previous contraction connecting
$C^*(K)$ with $H^*(K)$, it is possible to design an algorithm computing the
commutative $A_{\infty}$--algebra structure of $H^*(K)$ reflecting the complete rational
homotopy type of $K$ \cite{Kad98}. We will see in the
next section
that the homotopy equivalence data structure is also essential in computing cohomology operations.

\section{Steenrod Cohomology Operations}

Let us suppose $\Lambda={\bf Z}_p$ ($p$ being a prime),
then it is possible to construct an algebraic minimal model for
any finite simplicial complex $K$, in which the associated contraction $(f^*,g^*,\phi^*)$
connects $C^*(K)$ with its cohomology.
From this data
and the combinatorial formulae for Steenrod squares and
reduced $p$th powers \cite{Ste47,ES62} at cochain level in terms of
face operators  established in \cite{GR99,Gon00},
Steenrod cohomology operations can effectively be  computed.

For instance, the formula for
 the Steenrod reduced power

 $${\cal P}_1:H^*(X)\ra H^{*p-1}(X)$$
 at cochain level \cite{Gon00}   is:

\begin{eqnarray*}
P_1(c)(\sigma)=\sum_{j=1}^{p-1}\;\;\sum_{i=jq}^{(j+1)q-1}&&
(-1)^{(i+1)(q+1)+1}\\
&&\mu(c(v_0,\dots,v_q)\\
&&\ot   (v_q,\dots,v_{2q})\\
&&\vdots    \\
&&\ot   c(v_{(j-2)q},\dots,v_{(j-1)q})\\
&&\ot   c(v_{(j-1)q},\dots,v_{i-q},v_i,\dots,v_{(j+1)q-1})\\
&&\ot   c(v_{(j+1)q-1},\dots,v_{(j+2)q-1})\\
&&\vdots    \\
&&\ot   c(v_{(p-2)q-1},\dots,v_{(p-1)q-1})\\
&&\ot   c(v_{(p-1)q-1},\dots,v_{pq-1})\\
&&\ot   c(v_{i-q},\dots,v_i)\,)\\
\end{eqnarray*}
where $c$ is a $q$--cocycle, $\sigma=(v_0,v_1,\dots,v_{pq-1})$ is a $(pq-1)$--simplex such that
$v_0<v_1<\cdots<v_{pq-1}$ and $\mu$ is the product on ${\bf Z}_p$.
Therefore, for calculating the cohomology class ${\cal P}_1(\alpha)$ with
 $\alpha\in H^q(K)$,
we only have to compute $f^*P_1g^*(\alpha)$.

In the particular case of Steenrod squares,

$$Sq^i:H^*(K;{\bf Z})\rightarrow  H^{*+i}(K;{\bf Z}_2)\,,$$
we can express them in a matrix form due to the fact that these
cohomology operations are homomorphisms.
Moreover, the process of diagonalization of such matrices can give us  detailed
information about the kernel and image of these cohomology operations.

\section{Adem Secondary Cohomology Operations}

For attacking the computation of secondary cohomology operations, we will see in this section
that the homotopy operator $\phi^*$ of the contraction associated to an algebraic minimal model
of a simplicial complex $K$ is essential.

First of all, we shall indicate how Adem secondary cohomology operations

$$\Psi_q: N^q(K)\ra H^{q+3}(K;{\bf Z}_2)/Sq^2 H^{q+1}(K;{\bf Z})$$
can be constructed (see \cite{Ade52}).
$N^q(K)$ denotes the kernel of $Sq^2: H^q(K;{\bf Z})\ra H^{q+2}(K;{\bf Z}_2)$
These operations appear using the known relation:

 $$Sq^2Sq^2\alpha+Sq^3Sq^1\alpha=0$$
for any $\alpha\in H^*(K;{\bf Z})$. For this particular relation there exist cochain mappings

$$E_{j}:C^*(K\ti K\ti K\ti K)\ra C^{*-j}(K)$$ such that
mod $2$

$$(c\smile_{q-2}c)\smile_q  (c\smile_{q-2}c)
+(c\smile_{q-1}c)\smile_q  (c\smile_{q-1}c)=\delta E_{3q-3}(c^4)\,,$$
where $\smile_k$ is the cup--$k$ product \cite{Ste47}
and $c$ is a $q$--cochain.
Recall that, at cochain level,
$Sq^i(c)=c\smile_{j-i}c$ mod $2$,  where $c$ is a $j$--cocycle.
Then
$\Psi_q$
is defined at cochain level by

$$\psi_q(c)=b\smile_{i+1} b+b\smile_{i+2}\delta b+E_{3i+3}(c)
+\eta(c)\smile_{i-1}\eta(c)+\eta (c)\smile_i\delta \eta(c)\,,$$
where $c$ is a $q$--cocycle representative of a cohomology class  of
$N^q(K)$, $b$ is a $(q+1)$--cochain such that $c\smile_{q-2}c=\delta b$ and
$\eta(c)=\frac{1}{2}(c\smile_q c+c)$.

If ${\bf Z}_2$ is the ground ring,
formulae for computing cup--$i$ products are well--known
\cite{Ste47}.
A method for obtaining ``economical" formulae for $E_{3i+3}$
in terms of face operations is given in \cite{Gon00}. For example,

 \begin{eqnarray*}
E_3(c^4)(\sigma)
&=&\mu(c(v_0,v_2,v_3)\ot c(v_0,v_1,v_2)\ot c(v_3,v_4,v_5)\ot c(v_2,v_3,v_5)\\
    &&+c(v_0,v_4,v_5)\ot c(v_3,v_4,v_5)\ot c(v_0,v_1,v_2)\ot c(v_0,v_1,v_2)\\
    &&+c(v_0,v_1,v_5)\ot c(v_3,v_4,v_5)\ot c(v_1,v_2,v_3)\ot c(v_1,v_2,v_3)\\
    &&+c(v_0,v_1,v_2)\ot c(v_2,v_4,v_5)\ot c(v_2,v_3,v_4)\ot c(v_2,v_3,v_4)\\
    &&+c(v_0,v_1,v_2)\ot c(v_2,v_3,v_5)\ot c(v_3,v_4,v_5)\ot c(v_3,v_4,v_5))\,,
\end{eqnarray*}
where $c$ is a $2$--cochain,  $\sigma=(v_0,v_1,...,v_5)$ is a $5$--simplex such that $v_0<v_1<\cdots<v_5$
and $\mu$ is the product on ${\bf Z}_2$.
Therefore, the steps for computing $\Psi_q$ are the following:

\begin{itemize}
\item[1.] Take $\alpha\in N^q(K)$ making use of the diagonalization of the matrix of
$Sq^2$ in dimension $q$.
\item[2.] Compute $b=\phi^*Sq^2g^*(\alpha)$.
\item[3.] Compute $f^*\psi g^*(\alpha)$.
\end{itemize}
Note that it is very easy to prove that

$$g^*(\alpha)\smile_{q+2}g^*(\alpha)=\delta \phi^*Sq^2g^*(\alpha)\,,$$
using
the relation $1-g^*f^*= \phi^*\delta+\delta\phi^*$.

\section{Some Comments}

All these results can be given in a more general framework working with not necessarily finite
simplicial complexes. Nevertheless, a contraction of the
chain complex associated to the simplicial complex to its (co)homology must exist
in order to develop the method.

Concerning the complexity,
  obtaining a contraction of
a finite simplicial complex $K$ to its (co)homology can be done using
   Delfinado--Edelsbrunner incremental algorithm \cite{ELZ00} which
 runs in time as most cubic in the number of simplices of the complex if the group of
 coefficients is a field.
On the other hand, another datum to take into account  is the number of summands of
the formulae for computing
cohomology operations at cocycle level.  For example the number of summands of $P_1$ over a
$q$--cocycle $c$ and a $(pq-1)$-simplex $\sigma$
is $(p-1)q$.

Finally,
in order to obtain the image of
any cohomology operations at cochain level
 over a representative
cocycle using our formulae, we have to compute them over a basis of $C_*(K)$ in the desired dimension.
A way of decreasing the complexity of this
is to do a ``topological" thinning of the simplicial complex $K$
in order to obtain a thinned simplicial subcomplex $M_{\mbox{\scsz top}}(K)$ of $K$, such that
there exists a contraction of $C_*(K)$ to $C_*(M_{\mbox{\scsz top}}(K))$
For example, one way to construct it is using simplicial collapses \cite{For99}.
Then we can apply our machinery to compute cohomology operations in the thinned
simplicial complex $M_{\mbox{\scsz top}}(K)$ and the results can be
easily interpreted in the ``big" simplicial complex $K$.

\end{document}